\documentclass{article}

\usepackage{amssymb}
\usepackage{amsmath}
\usepackage{latexsym}
\usepackage{mathrsfs}
\usepackage{graphicx}
\usepackage[all]{xy}
\newtheorem{Th}{Theorem}

\newcommand{\map}[5]{#1 \quad : \quad #2\longrightarrow #3\qquad #4\longmapsto #5}

\title{Kernel-Cokernel Sequence for Composition and Its Applications}
\author{XIONG Rui}
\begin{document}
\maketitle

\begin{abstract}
  If this article, an elementary kernel-cokernel exact sequence is introduced for
  $A\stackrel{f}\to B\stackrel{g}\to  C$. Some relative sequences and applications are discussed.
  This result can simplify some proofs---the indices of Frodholm operators, Harada and Sai theorem, and the derived couples of exact couples.
\end{abstract}
%

\tableofcontents
\def\cok{\operatorname{cok}}
\def\im{\operatorname{im}}

\section{Sequence for composition}

\begin{Th}In an abelian category $\mathcal{C}$, if we have the morphism $A\stackrel{f}\to B\stackrel{g}\to  C$, then there is an exact sequence
$$0\to \ker f\to \ker gf\to \ker g\stackrel{\delta}\to \cok f\to \cok gf \to \cok g\to 0.$$
Here $\delta$ is the composition of $\ker g\to C\to \cok f$.
\end{Th}

I suggest to draw this diagram in the following way
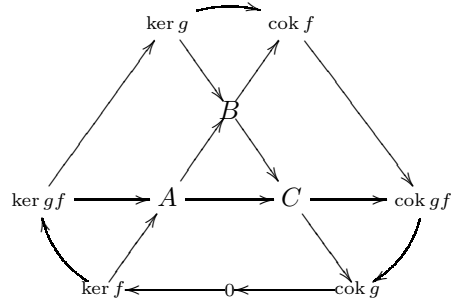
\def\B#1{{\displaystyle #1}}
\def\s#1{{\scriptstyle #1}}
$$\xymatrix{
\ar@{}[dr];[drr]|{\B{B}}="B"
\ar@{}[ddd];[dddr]|{\ker f}="kf"
\ar@{}[dddrr];[dddrrr]|{\cok g}="cg"
\ar@{}[dddr];[dddrr]|{0}="zero"
& \s{\ker g}\ar"B"\ar@/^/[r] & \s{\cok f} \ar[ddr]\\
& & \ar"B";[d] \ar"B";[u] \\
\s{\ker gf} \ar[r]\ar[uur] &A \ar[r]\ar"B" & C \ar[r]\ar"cg" & \s{\cok gf}\ar@/^/"cg" \\
\ar@/^/"kf";[u]\ar"kf";[ur]& \ar"zero";"kf"& \ar"cg";"zero"&
}$$

The proof is not difficult by diagram chasing. I'd like to give a proof using snake lemma as following
$$
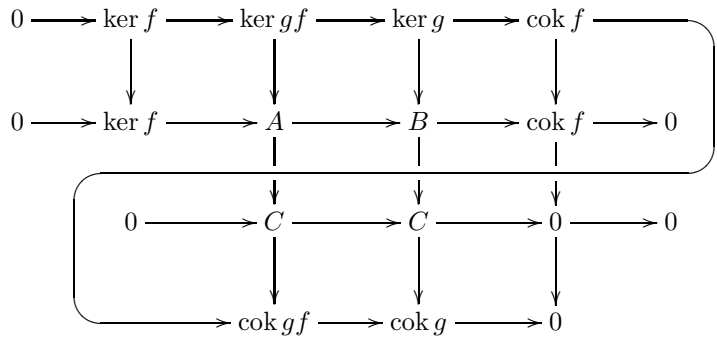
\xymatrix{
0\ar[r]&\ker f \ar[r]\ar[d] & \ker gf \ar[r]\ar[d] & \ker g \ar[r]\ar[d] & \cok f\ar[d]
\ar@{}[dr];[ddrr]|{}="a"
\ar@{}[ddlll];[dddllll]|{}="b"
\ar`"a"`"a"`"b"`[dddll][dddll]
& &\\
0\ar[r]&\ker f\ar[r] & A\ar[r]\ar[d]|{\hole}& B\ar[r]\ar[d]|{\hole} & \cok f\ar[d]|{\hole}\ar[r] & 0&\\
&0 \ar[r] &C\ar[d]\ar[r]& C\ar[d]\ar[r] & 0\ar[d]\ar[r] & 0&\\
&& \cok gf \ar[r]& \cok g \ar[r] & 0&
}$$

The exactness also follows from \emph{Wall theorem}, which is known more by topologists,
for example, \cite{bredon2013topology} Page 189.
$$
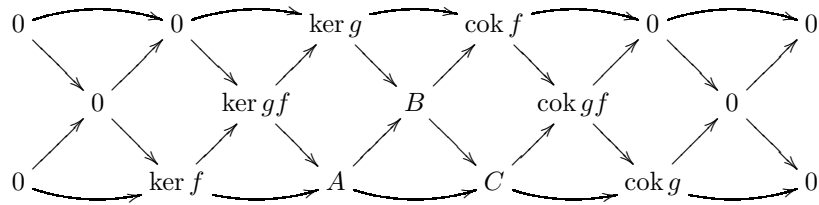
\xymatrix@!=0.5pc{
0 \ar@/^/[rr]\ar[dr]&&0\ar@/^/[rr]\ar[dr] && \ker g\ar@/^/[rr]\ar[dr] && \cok f\ar@/^/[rr]\ar[dr] && 0\ar@/^/[rr]\ar[dr] && 0\\
&0 \ar[ur]\ar[dr]&& \ker gf\ar[ur]\ar[dr] &&  B\ar[ur]\ar[dr] && \cok gf\ar[ur]\ar[dr] && 0\ar[ur]\ar[dr]\\
0 \ar@/_/[rr]\ar[ur]&&\ker f\ar@/_/[rr]\ar[ur]&& A\ar@/_/[rr] \ar[ur]&& C\ar@/_/[rr] \ar[ur] && \cok g\ar@/_/[rr]\ar[ur] &&0
}$$


Of course, above sequence covers some trivial algebraic conclusions
\begin{itemize}
  \item If $gf$ is injective, then $f$ is injective; if $gf$ is surjective, then $g$ is surjective.
  \item If $f$ and $g$ are injective (surjective), then so is $gf$.
\end{itemize}

For the diagram
$$
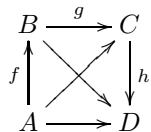
\xymatrix{
B\ar[r]^g\ar[dr]&C\ar[d]^h\\
A\ar[u]^f\ar[ur]\ar[r]&D
}$$
We have four kernel-cokernel sequences, it is not difficult to draw them in the following way.

\begin{Th}We have the following diagram
$$
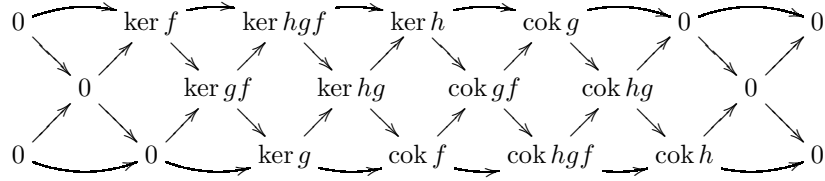
\xymatrix@!=0.1pc{
0 \ar@/^/[rr]\ar[dr]&&\ker f\ar@/^/[rr]\ar[dr] && \ker hgf\ar@/^/[rr]\ar[dr] && \ker h \ar@/^/[rr]\ar[dr] && \cok g\ar@/^/[rr]\ar[dr]
&& 0\ar@/^/[rr]\ar[dr] && 0\\
&0 \ar[ur]\ar[dr]&& \ker gf\ar[ur]\ar[dr] &&  \ker hg\ar[ur]\ar[dr] && \cok gf\ar[ur]\ar[dr] && \cok hg \ar[ur]\ar[dr]&&0\ar[ur]\ar[dr]\\
0 \ar@/_/[rr]\ar[ur]&&0\ar@/_/[rr]\ar[ur]&& \ker g\ar@/_/[rr] \ar[ur]&& \cok f\ar@/_/[rr] \ar[ur] && \cok hgf\ar@/_/[rr]\ar[ur]
&&\cok h \ar@/_/[rr]\ar[ur]&& 0
}$$
\end{Th}

\section{Sequence for commutative square}

For
$$
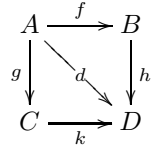
\xymatrix{
A\ar[r]^f\ar[d]_g\ar[dr]|{d} &B \ar[d]^h\\
C\ar[r]_k &D}$$
We have two kernel-cokernel sequences, the problem is how they intersect?

\begin{Th}We have the following diagram
$$
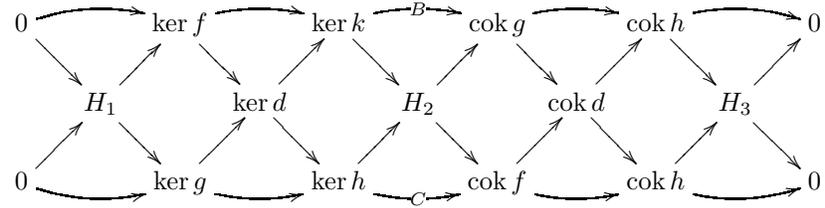
\xymatrix@!=0.5pc{
0\ar@/^/[rr]\ar[dr]&&\ker f \ar@/^/[rr]\ar[dr] && \ker k\ar@/^/[rr]|{B}\ar[dr] && \cok g \ar@/^/[rr]\ar[dr] && \cok h\ar@/^/[rr]\ar[dr] && 0\\
&H_1 \ar[ur]\ar[dr]&& \ker d \ar[ur]\ar[dr] &&  H_2\ar[ur]\ar[dr] && \cok d\ar[ur]\ar[dr] && H_3\ar[ur]\ar[dr]\\
0\ar@/_/[rr]\ar[ur]&&\ker g \ar@/_/[rr]\ar[ur]&& \ker h\ar@/_/[rr]|{C} \ar[ur]&& \cok f\ar@/_/[rr] \ar[ur] && \cok h\ar@/_/[rr]\ar[ur] && 0
}$$
with each braid exact. Where $H_1,H_2,H_3$ are exact the homology groups of
$$0\to \stackrel{1}A\to \stackrel{2}{B\oplus C} \to \stackrel{3}D\to 0. $$
\end{Th}

The proof is directly by diagram chasing.

In some special cases, the above exact net is very useful to prove the some results diagram-chasing-consuming.
\begin{itemize}
  \item When the square is pull-back, then $H_1=H_2=0$, so $\ker f=\ker k$ and $\cok f\to \cok k$ is surjective.
  \item Conversely, when the natural map $\ker f\to \ker k$ is isomorphism, and the natural map $\cok g\to \cok h$ is surjective, the square is pull-back.
  \item When the square is push-out, then $H_2=H_3=0$, so $\cok f=\cok k$ and $\ker f\to \ker k$ is injective.
  \item Conversely, when the natural map $\cok f\to \cok k$ is isomorphism, and the natural map $\ker g\to \ker h$ is injective, the square is push-out.
  \item When $k$ is injective, $H_1=\ker f$. Furthermore, $H_2$ and $H_3$ are exactly the kernel and cokernel of $\cok f\to \cok k$.
  Now the long exact sequence
  $$0\to H_1 \to \ker g\to \ker h\to H_2 \to \cok g\to \cok h\to H_3\to 0$$
  is exactly the sequence snake lemma claims for the following diagram
  $$\xymatrix{
  \ker f\ar[r]& A\ar[r]^f\ar[d]_g &B \ar[d]^h\ar[r] & \cok f \ar[d]\ar[r] & 0\\
  0\ar[r] & C\ar[r]_k &D \ar[r]& \cok k \ar[r] & 0
  }$$
  \item When $f$ is surjective, $H_3=\cok h$. Furthermore, $H_1$ and $H_2$ are exactly the kernel and cokernel of $\ker f \to \ker k$.
  Now the long exact sequence
  $$0\to H_1 \to \ker g\to \ker h\to H_2 \to \cok g\to \cok h\to H_3\to 0$$
  is exactly the sequence snake lemma claims for the following diagram
  $$\xymatrix{
  0\ar[r] & \ker f\ar[r]\ar[d]& A\ar[r]^f\ar[d]_g &B \ar[d]^h\ar[r] &  0\\
  0\ar[r] & \ker k\ar[r] & C\ar[r]_k &D \ar[r]& \cok k
  }$$
\end{itemize}

\section{Index of Fredholm operator or Herbrand quotient}

\def\ind{\operatorname{ind}}

Let $\mathcal{C}$ be an abelian category, if there is a well-defined length function $\ell$
over some Serre subcategory (closed under kernel and cokernel),
say $\mathcal{C}_{\ell<\infty}$ such that for any exact sequence in $\mathcal{C}_{\ell<\infty}$,
$$0\to A'\to A\to A''\to 0.$$
we have $\ell(A')+\ell(A'')=\ell(A)$.
Or abstractly, $\ell$ is an element in Grothendieck group of $\mathcal{C}_{\ell<\infty}$.
We say an morphism $A\stackrel{f}\to B$ is Fredholm, if $\ker f$ and $\cok f$ are in $\mathcal{C}_{\ell<\infty}$.

For any Fredholm morphism $A\stackrel{f}\to B$, we can consider the index (Herbrand quotient)
$$\ind (f)=\ell(\cok f)-\ell(\ker f). $$

Using our exact sequence, and the euler characteristic we have
\begin{Th}If $f$, $g$ and $gf$ are all Fredholm morphisms, then
$$\ind (gf)=\ind f+\ind g. $$
\end{Th}

Note that the proof that the index of composition of Fredholm operators is the sum of indices is really a hard work in standard books.
In functional analysis, we define the Fredholm operator to be the bounded linear operator between Banach space with closed image and
having finite dimensional kernel and cokernel. the index of Fredholm operator $A$ is defined to be the
$$\ind A=\dim \ker A-\dim \ker A^*, $$
which is exactly $\dim \ker A-\dim \cok A$, see \cite{conway2010course} Page 354.
But the proof of $\ind (BA)=\ind A+\ind B$ there is really hard.
When we proved the composition of Fredholm operator is still Fredholm operator (this part is easy),
the proof can be direct by our more general consequence.



\bigbreak
It also works well for the (co)homology of finite groups,
since the (co)homology groups are finite for finite generated module.
Here it is efficient to take $h(f)=\frac{\#(\cok f)}{\#(\ker f)}$ rather than additive.
Let $M$ be an $G$-module, consider the norm map
$$\map{N_G^M}{M_G}{M^G}{x}{\sum gx}$$
we can define
$$h(M)=h(N_G)=\frac{\# \tilde{H}^{0}(G;M)}{\# \tilde{H}_{-1}(G,M)}. $$
Here $\tilde{H}^\bullet(G;M)$ is Tate-cohomology group.
In particular, it works very well for cyclic groups, since now Tate group is cyclic, now we have $h(M)=h(M')+h(M'')$ for exact sequence
$$0\to M'\to M\to M''\to 0. $$
See \cite{serre1979local} Page 133.
%
%

\section{A simpler proof of Harada and Sai theorem}

In representation of associative algebras,
we have the following theorem by Harada and Sai, which is useful in proving the first Brauer–Thrall conjecture
c.f. \cite{simson2007elements} volume one page 139.
The proof there is standard, not difficult, but I will give a direct and simpler proof.

\begin{Th}[Harada and Sai]Let $M_1,\ldots,M_{2^n}$ be indecomposable module of length no more than $n$,
if we have a sequence
$$M_1\stackrel{f_1}\to M_2\to \cdots \to M_{2^n-1}\stackrel{f_{2^n-1}}\to M_{2^n}$$
if each of $f_i$ is not isomorphism, then the composition $f_{{2^n}-1}\cdots f_1=0$.
\end{Th}

We prove by induction that the image of $f_{{2^m}-1}\cdots f_1=0$ them are of length no more than $n-m$.
For $m=1$, this is trivial, since it is not an isomorphism.
For general case, consider the half $M_1\stackrel{g}\to M_{2^{m-1}}\stackrel{h}\to M_{2^m}$. We have
$$0\to \ker g\stackrel{*}\to  \ker hg\to \ker h\stackrel{\delta}\to \cok g\to \cok hg\stackrel{**}\to \cok h\to 0$$
\begin{itemize}
\item If $(*)$ is not isomorphism, then $\im hg=M_1/\ker hg$ whose length is strictly less than $\im g=M_1/\ker g$, so it follows by induction hypothesis.
\item If $(**)$ is not isomorphism, then $\im hg$ whose length is strictly less than $\im h$, so it follows by induction hypothesis.
\item If $\delta$ is isomorphism, then $\ker h\to M_{2^{m-1}}\to \im g$ is an isomorphism, but $M_{2^{m-1}}$ is assumed to be indecomposable,
so $\ker h=0$ and $\im g=0$. Now $M_{2^{m-1}}$ is of less length of $M_{2^n}$---otherwise they are isomorphic.
\end{itemize}

\section{Cubic zero and Quartic zero sequences}

The purpose of this section is to give a natural way to find why the derived couple of exact couple is exact.

\bigbreak
We say a sequence
$$(ABC)_\bullet=\cdots \to A_1\stackrel{d}\to B_1 \stackrel{d}\to C_1 \stackrel{d}\to A_2\stackrel{d}\to B_2\to \cdots$$
is cubic zero if it satisfies $d^3=0$.

It is more efficient to write
$$
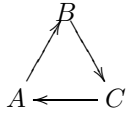
\xymatrix{
\ar@{}[r]|{\B{B}}="B"&\ar"B";[d]\\
A\ar"B"& C \ar[l]
}$$
Denote $H^X_Y$ the homology group at $X$ in the complex made up by $X$ and $Y$ for $X,Y\in \{A,B,C\}$.

\begin{Th}
We have the following long exact sequence
$$
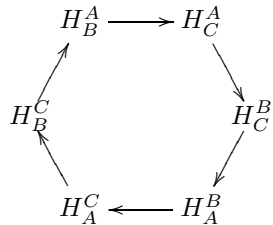
\xymatrix{
\ar@{}[d];[dr]|{\B{H^C_{B}}}="L"
\ar@{}[drr];[drrr]|{\B{H^B_{C}}}="R"
& H^A_{B} \ar[r] & H^A_{C}\ar"R"\\
\ar"L";[ur]&&& \ar"R";[dl]\\
&H^C_{A}\ar"L"&H^B_{A}\ar[l]\\
}$$
\end{Th}

The proof is of course can be done by diagram chasing,
but the abstract proof is not difficult, parallel to the proof of long exact sequence for short exact sequence.
$$\xymatrix{
&& H^B_A \ar[r]\ar[d] & H^C_A \ar[r]\ar[d] & H^B_C\ar[d]
\ar@{}[dr];[ddrr]|{}="a"
\ar@{}[ddlll];[dddllll]|{}="b"
\ar`"a"`"a"`"b"`[dddll][dddll] &\\
&& \cok[A\to B]\ar[r]\ar[d]|{\hole} & \cok[A\to C]\ar[r]\ar[d]|{\hole} &  \cok[C\to B]\ar[d]|{\hole} \ar[r]& 0&\\
&0\ar[r] & \ker[A\to B]\ar[r]\ar[d] & \ker[A\to C]\ar[r]\ar[d] &  \ker[C\to B]\ar[d]&&\\
& & H^A_B \ar[r] & H^A_C\ar[r] & H^C_B
}$$

Consider we say a sequence
$$(ABCD)_\bullet=\cdots \to A_1\stackrel{d}\to B_1 \stackrel{d}\to C_1 \stackrel{d}\to D_1 \stackrel{d} \to A_2\stackrel{d}\to B_2\to \cdots$$
is quartic zero if it satisfies $d^4=0$.
Still we denote $H^X_Y$ the homology group at $X$ in the complex made up by $X$ and $Y$ for $X,Y\in \{A,B,C,D\}$.

\begin{Th}[Quartic zero Lemma]
We have the following commutative diagram
$$\xymatrix@!=0pc{
\cdots \ar@/^/[rr]\ar[dr]&&H_B^A\ar@/^/[rr]\ar[dr] && H^A_D\ar@/^/[rr]\ar[dr] && H^C_D\ar@/^/[rr]\ar[dr] && H^C_B\ar@/^/[rr]\ar[dr] && \cdots\\
& H_B^D\ar[ur]\ar[dr]&& H_C^A\ar[ur]\ar[dr] &&  H_D^B\ar[ur]\ar[dr] && H_A^C\ar[ur]\ar[dr] && H_B^D\ar[ur]\ar[dr]\\
\cdots \ar@/_/[rr]\ar[ur]&&H_C^D\ar@/_/[rr]\ar[ur]&& H_C^B\ar@/_/[rr] \ar[ur]&& H_A^B\ar@/_/[rr] \ar[ur] && H^D_A\ar@/_/[rr]\ar[ur] &&\cdots
}$$
\end{Th}

It is efficient to draw it as octahedron
$$\xymatrix@!=2pc{
\ar@{}[dr]|!{[dd];[r]}{\B{BC}}="U"
\ar@{}[ddrr]|!{[ddr];[rr]}{\B{AD}}="D"
& AB\ar[dr]\ar"D" & \\
BD\ar[dr]\ar[ur] && AC \ar"U"\ar"U";[ul]\ar"U";[ll] \ar"D"\\
& CD \ar[ur]\ar"U" \ar"D";[ul]\ar"D";[]&
}$$
%

\begin{Th}[Quadratic zero Lemma]
Let $C$ be a differential object, assume $d$ is decomposed into $C\stackrel{e}\to D\stackrel{f}\to C$, 
with $D\stackrel{f}\to C\stackrel{e}\to D$ exact. 
Then we have the following long two exact sequences
$$\begin{array}{c}
\cdots \to H(C)\to \ker f\to \cok e \to H(C)\to \cdots\\
\cdots \to H(C)\to \cok f\to \ker e \to H(C)\to \cdots\\
\end{array}$$
\end{Th}

Since firstly, we have the following diagram
$$\xymatrix@!=0pc{
\cdots \ar@/^/[rr]|{0}\ar[dr]&&\ker f\ar@/^/[rr]\ar[dr] &&\cok e\ar@/^/[rr]|{0}\ar[dr] && \ker f\ar@/^/[rr]\ar[dr] && \cok e\ar@/^/[rr]|{0}\ar[dr] && \cdots\\
& H(C)\ar[ur]\ar[dr]&& D\ar[ur]\ar[dr] &&  H(C)\ar[ur]\ar[dr] && D\ar[ur]\ar[dr] && H(C)\ar[ur]\ar[dr]\\
\cdots \ar@/_/[rr]|{0}\ar[ur]&&\cok f\ar@/_/[rr]\ar[ur]&& \ker e\ar@/_/[rr]|{0} \ar[ur]&& \cok f\ar@/_/[rr] \ar[ur] && \ker e\ar@/_/[rr]|{0}\ar[ur] &&\cdots
}$$
The upper one is zero since it is induced by $ef$, the lower one is zero by exactness. 
Then by direct diagram chasing, we get the desired long exact sequence. 

As some readers like, we have the following commutative diagram
$$\xymatrix@!=0.2pc{
\cdots & \cok f \ar[dddr]\ar[rr]&& \ker e \ar[dr] &\to\!\! H\!\! \to
& \cok f \ar[dddr]\ar[rr]&& \ker e \ar[dr]& \to\!\! H\!\! \to
& \cok f \ar[dddr]\ar[rr]&& \ker e \ar[dr] & \cdots \\
C\ar[ur]\ar[ddrr]\ar[rrrr] &&&& 
C\ar[ur]\ar[ddrr]\ar[rrrr] &&&& 
C\ar[ur]\ar[ddrr]\ar[rrrr] &&&& C\\ \\
&& D\ar[uurr]\ar[uuur]\ar[dr] && 
&& D\ar[uurr]\ar[uuur]\ar[dr] && 
&& D\ar[uurr]\ar[uuur]\ar[dr] && \\
\cdots& \ker f\ar[ur]\ar[rr] && \cok e & \to \!\!H\!\! \to
& \ker f\ar[ur]\ar[rr] && \cok e & \to\!\! H\!\! \to
& \ker f\ar[ur]\ar[rr] && \cok e & \cdots
}$$
%

\section{Derived couple of Exact couple}

Now, let us give a ``natural'' proof, that the derived couple of exact couple is exact, c.f. \cite{weibel1995introduction} page 153.
The proof can be directly from elements by elements check, but there is no evidence why it is true.
But as we established in the previous section, it naturally appears.

Assume we have an exact couple
$$\xymatrix{
\ar@{}[d];[dr]|{\B{E}}="E"
D\ar[r]^{\alpha} & D \ar"E"^{\beta}\\
\ar"E";[u]^{\gamma}&&
}$$
It defines an differential object $E$ with $\delta=\beta\gamma$.
It satisfies the assumption above, so we have a long exact sequence,
$$\cdots\to H(C)\to \ker \beta\to \cok \gamma \to H(C)\to\cdots.$$
Since our assumption that it is an exact couple, $\ker \beta=\im \alpha=\cok \gamma$.
We have a new exact couple now
$$\xymatrix{
\ar@{}[d];[dr]|{\B{H(E)}}="E"
\im \alpha\ar[r]^{\alpha} & \im \alpha \ar"E"^{\beta}\\
\ar"E";[u]^{\gamma}&&
}$$

\section{Remarks}

In someway, the cubic zero sequence introduced in the article can be viewed as a analogy of mapping cone,
and it satisfies very trivial rotation axiom, it may form a triangulated structure in someway,
the quartic zero Lemma can be thought as octahedral axiom.

I found this sequence when I was a third year student in university in 2018.
But I did not realize that it is useful to give simplified proofs as we present in above sections.
The original proof (written in standard text book) is not easy and not direct.
I thought a lot of mathematicians do not know this easy conclusion.
I think it may be useful to diffuse this conclusion.
\vfill
\bibliographystyle{plain}
\bibliography{bibfile}

\begin{thebibliography}{1}

\bibitem{bredon2013topology}
Glen~E Bredon.
\newblock {\em Topology and geometry}, volume 139.
\newblock Springer Science \& Business Media, 2013.

\bibitem{conway2010course}
John~B Conway.
\newblock {\em A course in functional analysis}, volume~96.
\newblock Springer, 2010.

\bibitem{serre1979local}
Jean-Pierre Serre.
\newblock Local class field theory.
\newblock In {\em Local Fields}, pages 188--203. Springer, 1979.

\bibitem{simson2007elements}
Daniel Simson and Andrzej Skowronski.
\newblock {\em Elements of the representation theory of associative algebras}.
\newblock Cambridge Univ. Press, 2007.

\bibitem{weibel1995introduction}
Charles~A Weibel.
\newblock {\em An introduction to homological algebra}.
\newblock Number~38. Cambridge university press, 1995.

\end{thebibliography}

\end{document}